\documentclass[preprint, 12pt]{elsarticle}
\usepackage[T1]{fontenc}
\usepackage{amsmath}
\usepackage{amsfonts}
\usepackage{amssymb}
\usepackage{amsthm}
\usepackage{eucal}
\usepackage{mathrsfs}
\usepackage{geometry}
\usepackage{hyperref}
\usepackage{theoremref}
\usepackage[dvipsnames]{xcolor}
\usepackage{float}
\usepackage{caption}
\usepackage{subcaption}

\usepackage{etoolbox}
\usepackage{multicol}
\usepackage{color}
\usepackage{xspace}
\usepackage{enumerate}
\makeatletter
\numberwithin{equation}{section}
\usepackage{setspace}
\DeclareUnicodeCharacter{0301}{\'{e}}
\allowdisplaybreaks
\usepackage[pagewise]{lineno}
\begin{document}
\title{\textbf{Construction and Box-counting Dimension of the Edelstein Hidden Variable Fractal Interpolation Function}}
\author[1]{Aiswarya T\corref{cor1}%
\fnref{fn1}}
\ead{aiswaryasidhu8113@gmail.com}
\author[1]{Srijanani Anurag Prasad\fnref{fn1}}
\ead{srijanani@iittp.ac.in}
\cortext[cor1]{Corresponding author}
\affiliation[1]{organization={Department of Mathematics and Statistics, Indian Institute of Technology Tirupati},
%addressline={Jawahar Nagar},
%city={Trivandrum},
% citysep={}, % Uncomment if no comma needed
% % between city and postcode
%postcode={695013},
%state={Kerala}
country={India}}
\date{}
\fntext[fn1]{Department of Mathematics and Statistics, IIT Tirupati, Yerpedu P.O., India, 517619}

\newtheorem{proposition}{Proposition}
\newtheorem{theorem}{Theorem}
\newtheorem{definition}{Definition}
\newtheorem{example}{Example}
\newtheorem{remark}{Remark}
\newtheorem{corollary}{Corollary}
\begin{abstract}
This paper presents the construction of a hidden variable fractal interpolation function using Edelstein contractions in an iterated function system based on a finite collection of data points. The approach incorporates an iterated function system where variable functions act as vertical scaling factors leading to a generalised vector-valued fractal interpolation function. Furthermore, the paper rigorously examines the smoothness of the constructed function and establishes an upper bound for the box-counting dimension of its graph.
\end{abstract}
\begin{keyword}
    Iterated Function Systems \sep Edelstein Contraction \sep Fractal Interpolation Function \sep Attractor \sep Box Dimension
    \hspace{2mm}
\end{keyword}
\maketitle
\section{Introduction}
The fractal interpolation method is a relatively new approach in approximation theory offering significant advantages for modelling datasets. Unlike traditional interpolation methods, fractal interpolation generates approximating functions that inherently capture nonlinearity and self-similarity.

Barnsley developed the notion of Fractal Interpolation Function(FIF) \cite{barnsley1986fractal} based on the idea of a function whose graph coincides with the attractor of an Iterated Function System(IFS) of Banach contractions. Subsequent studies have expanded the idea of FIFs. These advancements led to the development of generalised FIFs, further enhancing the flexibility and applicability of the method across complex datasets. This includes the extension to multivariate cases \cite{dalla2002bivariate}, \cite{xie1997study}, %\cite{bouboulis2007fractal} 
and the development of Super FIFs\cite{srijanani15}. Another extension of FIFs incorporated generalised contractions such as Rakotch and Matkowski contractions\cite{ri2017new},\cite{ri2019new} within the iterated function systems to construct nonlinear FIFs. The most generalised construction form to date is the construction of a nonlinear FIF with Edelstein contractions in the underlying iterated function system\cite{pasupathi2024very}. A related development is the construction of histopolating fractal functions, which adapt the fractal approach to approximate histogram-type data instead of pointwise interpolation \cite{barnsley2023histopolating}.

In parallel, the Hidden Variable Fractal Interpolation Function(HVFIF) was introduced as a function whose graph represents the projection of the attractor in a subspace of a particular IFS\cite{barnsley1989hidden}. This construction offered greater flexibility for the fractal functions. This led to the construction of functions such as coalescence hidden variable FIFs\cite{chand2005coalescence}, hidden variable recurrent FIFs\cite{yun2019hidden} and hidden variable FIFs based on iterated function systems consisting of Rakotch contractions\cite{kim2019construction}. Various properties of such functions - including dimension\cite{barnsley1989hidden}, smoothness %\cite{yun2020box},
\cite{ri2021smoothness},\cite{kapoor2009smoothness}, stability\cite{yun2020box},\cite{kapoor2009stability},%\cite{kim2019construction} 
and monotonicity properties\cite{katiyar2016hidden} - were studied. The construction of parametrized $\mathcal{A}-$ fractal functions based on HVFIF and approximation properties are also discussed in~\cite{chand2015approximation}.

The article is divided into the following sections: Section~\ref{prelims}  presents the definitions of some generalised contractions and fixed point theorems related to these contractions. This section includes the definition of an iterated function system and some results regarding iterated function systems with generalised contractions. Section~\ref{constrn} discusses the construction of the Edelstein hidden variable fractal interpolation function and includes examples to demonstrate the process. Section~\ref{smoothness} examines the smoothness of the constructed function. Section~\ref{boxdimn} provides an upper bound for the box-counting dimension.
\section{Preliminaries}\label{prelims}
In this section, definitions of some generalised contractions are presented, along with fixed-point theorems concerning these contractions. Furthermore, a few results related to iterated function systems consisting of generalised contractions are also discussed.

Let $(X,\mu)$ be a metric space. The Banach contraction is one of the most commonly discussed contractions.
\begin{definition}
\cite{barnsley2014fractals} Let $h$ be a function on $X$. Suppose there is a $k\in[0,1)$ such that
\begin{equation}
    \nonumber \mu(h(x),h(y))\leq k\mu(x,y), \forall \ x,y\in X.
\end{equation}
Then $h$ is a Banach contraction.
\end{definition}
The fixed point theorem for Banach contractions is considered a fundamental result in functional analysis.
\begin{proposition}
    \cite{banach1922operations} If $g$ is a Banach contraction on a complete metric space $X$, then there exists a unique fixed point for $g$ in $X$. 
\end{proposition}
Numerous classes of contractions extend the concept of Banach contractions and the Banach fixed-point theorem. One such class is the Edelstein contractions. 
\begin{definition}
\cite{Edelstein1962} If a function $g$ on $X$ satisfies
\begin{equation*}
    \mu(g(x),g(y))<\mu(x,y),~\forall \  x,y\in X,\text{ with }x\neq y,
\end{equation*}
then $g$ is an Edelstein contraction.
\end{definition}
It is established that a fixed point exists for an Edelstein contraction within the framework of a compact metric space. 
\begin{proposition}\cite{berinde2007iterative}
    Every Edelstein contraction $g$ on a compact metric space $X$ has a unique point $x'\in X$ such that $g(x')=x'$.
\end{proposition}
The construction of any FIF relies on an iterated function system and its attractor.
\begin{definition}    \cite{barnsley2014fractals} An iterated function system, $\mathcal{I}$, is a complete metric space $X$ with a finite number of contraction maps, $f_1, f_2,\ldots,f_N,$ defined on $X$ and is represented as $\mathcal{I}=\left\{X;f_j,j=1,2,\ldots,N\right\}$.
\end{definition}
 Let $\mathcal{H}(X)$ denote the set of all non-empty compact subsets of the metric space $X$.
\begin{definition}
    \cite{barnsley2014fractals} A set $A\in \mathcal{H}(X)$ is an attractor for IFS ~$\mathcal{I}$ if it is the unique set in $\mathcal{H}(X)$ satisfying
    \begin{equation*}
        A=\bigcup_{j=1}^N f_j(A).
    \end{equation*}
\end{definition}
\begin{proposition}\label{edifs}\cite{pasupathi2024very}
Let $(X,\mu)$ be a compact metric space and $\mathcal{I}$ be an IFS with Edelstein contractions. Then, $\mathcal{I}$ possess an attractor.
\end{proposition}
Box-counting dimension is one of the most widely used dimensions.
\begin{definition}\cite{falconer2013fractal}
    Let $F$ be a bounded subset of $\mathbb{R}^n$ and $N_{\epsilon}$ be the number of $\ epsilon$-mesh
squares that intersect $F$. Then, the upper box-counting dimension of $F$ is defined as
$$\overline{dim}_B(F)=\overline{\lim\limits_{\epsilon\rightarrow 0}}\frac{log N_\epsilon}{-log \epsilon}.$$
\end{definition}
\section{Construction of Edelstein Hidden Variable Fractal Interpolation Function}\label{constrn}
The hidden variable fractal interpolation function based on Edelstein contractions is constructed in this section.Let $V= \left\{(t_j,v_j)\in I\times \mathbb{R} : j=0,1,\ldots N \right\}$ be the given data set where~ $t_0<t_1<\ldots <t_N$. The generalised data set corresponding to $V$ is
\begin{equation*}
    V^*=\left\{(t_j,v_j,w_j)\in I\times \mathbb{R}^2 : j=0,1,\ldots N \right\},
\end{equation*}
Denote $[t_0,t_N]$ by $I$ and  $[t_{j-1}, t_j]$ as $I_j$ for all $j=1,2,\ldots N$. Let $L_j:I\rightarrow I_j$ be contractive homeomorphisms given by $L_j(t)=a_jt+f_j$ such that $L_j(t_0)=t_{j-1}$ and $L_j(t_N)=t_j$ for all $j\in \left\{1,2,\ldots, N\right\}$.
For any real valued function $g$, denote $\|g\|_{\infty}=\sup\limits_t|g(t)|$ and $l_g$ be the Lipschitz constant for a Lipschitz function $g$. Let $b_j,c_j,d_j$ and $e_j$ be real valued Lipschitz functions on $I$ such that 
\begin{equation}
    \|b_j\|_{\infty}+\|d_j\|_{\infty}\leq 1 \text{ and } \|c_j\|_{\infty}+\|e_j\|_{\infty}\leq 1\ \forall \ j\in\left\{1,2,\ldots N\right\}.
\end{equation}
Let $s_j$ and $r_j$ be functions on $\mathbb{R}$ such that they are Edelstein contractions. Further, define $F_j:I\times \mathbb{R}^2 \rightarrow \mathbb{R}^2$ as
\begin{equation}
    F_j(t,v,w)=D_j(t)S_j(v,w)+Q_j(t) \label{fidefn}
\end{equation}
where,
\begin{equation}
 D_j(t)= \begin{bmatrix}
b_j(t) & c_j(t)\\
d_j(t) & e_j(t)
\end{bmatrix}, \quad 
S_j(v,w)= \begin{bmatrix}
s_j(v) \\
r_j(w) 
\end{bmatrix}, \ \text{ and } 
Q_j(t)= \begin{bmatrix}
p_j(t) \\
q_j(t) 
\end{bmatrix}~\text{ for }~j=1,\ldots N.
\end{equation}\label{fidefnexpn}
Here, $p_j(t)$ and $q_j(t)$  are real-valued Lipschitz functions defined on $I$ such that
\begin{equation*}
    F_j(t_0,v_0,w_0)=(v_{j-1},w_{j-1}) \text{ and } F_j(t_N,v_N,w_N)=(v_{j},w_{j})~\text{ for }~j=1,2,\ldots ,N.
\end{equation*}
$F_j$ are Edelstein contractions with respect to the Manhattan metric in the second and third variables.Let $K\subset \mathbb{R}^2$ be the compact set containing $\left\{(v_j,w_j)\ :\ j=0,1,\ldots N\right\}$ such that $F_j(t,v,w)\in K$ for all $(t,v,w)\in I\times K$. Define the IFS $\mathcal{I}=\left\{I\times K, g_j, j=1,\ldots N\right\}$ where, 
\begin{equation}
    g_j(t,v,w)=(L_j(t),F_j(t,v,w))~\text{ for all }~(t,v,w)\in I\times K~\mbox{and} \ j=1,2,\ldots, N.\label{defw_i}
\end{equation}
Consider the space of $\mathbb{R}^2$ valued functions on $I$, $\mathcal{C}(I)$ equipped with the metric $d_C$ induced by the Manhattan metric.Let $C_e(I)=\{h\in C(I):h(t_0)~=~(v_0,w_0),~ h(t_N)~=~(v_N,w_N) \}$ and $C_d(I)=\{ h\in C_e(I)~:~h(t_j)~=~(v_j,w_j) ,\\ j =1,2,\ldots N-1 \}$.
Then $C(I)$, $C_e(I)$ and $C_d(I)$ are complete metric spaces with respect to the metric~$d_C$. Define the Read-Bajraktrevi\'{c} (R.B.) operator $\mathcal{R} : C(I)\rightarrow C(I)$ as
\begin{equation}
    \mathcal{R}h(t)=F_j(L_j^{-1}(t),h(L_j^{-1}(t)))~\text{ for }~t\in I_j~\text{ and }~j=1,2,\ldots, N. \label{RBdefn}
\end{equation}
\begin{theorem}\label{thrbfp}
Let $h$ be a function in $C_e(I)$. Then $\mathcal{R}h$ belong to $C_d(I)$. Moreover, $\mathcal{R}$ has a fixed point $f$ in $C_d(I)$.
\end{theorem}
\begin{theorem}
    The maps $g_j$ defined in Eq.~\eqref{defw_i} are Edelstein contractions with respect to the metric equivalend to the Euclidean metric. Moreover, the graph of the fixed point of the R.B. operator $\mathcal{R}$ defined in Eq.~\eqref{RBdefn} is the attractor of the IFS $\mathcal{I}$.
    \end{theorem}
    \begin{definition}[Edelstein Hidden Variable Fractal Interpolation Function]
    The first component of the vector-valued continuous function $f$ denoted as $f_1$ is the Edelstein hidden variable fractal interpolation function(EHVFIF) associated with the data set $\left\{(t_j, v_j)\in \mathbb{R}^2: j=0,1,\ldots N \right\}$.
\end{definition}
\section{Smoothness of the Edelstein Hidden Variable Fractal Interpolation Function}\label{smoothness}
\begin{theorem}\label{f_1smoothness}
 Consider the EHVFIF $f_1:I\rightarrow\mathbb{R}$ as defined in the preceding section. For appropriate choices of entries in $D_j,~\forall j=1,\ldots,N$ there exists positive real numbers $K$ and $\alpha$ such that
  \begin{equation}
     |f_1(t)-f_1(t')|<K|t-t'|^{\alpha}\quad \forall \ t,t'\in I.\label{holderf_1}
  \end{equation}
 \end{theorem}
 \section{Box Dimension of the Edelstein Hidden Variable Fractal Interpolation Function}\label{boxdimn}
 \begin{theorem}
Consider the EHVFIF $f_1: I \rightarrow \mathbb{R}$ as defined in Theorem~\ref{f_1smoothness}. Let $\alpha$ denote the same quantity as obtained in Theorem~\ref{f_1smoothness}. Then 
$$\overline{dim}_B(G(f_1))\leq 2-\alpha.$$
\end{theorem}
\bibliographystyle{unsrt}
\bibliography{EHVFIF}
\end{document}